\documentclass[12pt]{article}
\usepackage[mathscr]{eucal}
\usepackage{amsbsy}
\usepackage{amssymb}
\usepackage{amsmath}
\usepackage{amsthm}
\usepackage{graphics}
\usepackage{epsfig}

\usepackage{latexsym}

\usepackage[left=2cm, top=2cm, right=1.5cm, bottom=2cm] {geometry}

\def\phi{\varphi}

\def\bbr{{\mathbb R}}

\def\bz{{\bf Z}}

\def\bz{{\mathbb Z}} 

\def\Bl1{{\bf 1}}
\def\B2{{\bf 2}}
\def\B0{{\bf 0}}

\def\a{\alpha}
\def\b{\beta}
\def\d{\delta}

\def\e{\varepsilon}
\def\g{\gamma}

\def\l{\lambda}

\def\=A8{\"o}

\def\fdd{\stackrel{f.d.d.}{\longrightarrow}}

\newcommand{\PP}{\mathbf{P}}

\newcommand{\beq}{\begin{equation}}
\newcommand{\eeq}{\end{equation}}
\newcommand\beqn{\begin{displaymath}}  
\newcommand\eeqn{\end{displaymath}}

 \newcommand{\NN}{\mathbb{N} }

\newcommand{\ind}[1]{\mathbbm{1}_{#1}}

\newcommand{\halmos}{\vspace{3mm} \hfill \mbox{$\Box$}\\[2mm]}
\theoremstyle{plain}
\newtheorem{teo}{Theorem}

\newtheorem{prop}[teo]{Proposition}
\theoremstyle{definition}
\newtheorem{definition}[teo]{Definition}

\newtheorem{remark}[teo]{Remark}

\usepackage{bbm}
\begin{document}

\title{Limit theorems for linear processes with tapered innovations and filters
\footnotemark[0]\footnotetext[0]{ \textit{Short title:}
linear processes with tapered innovations }
\footnotemark[0]\footnotetext[0]{%
\textit{MSC 2010 subject classifications}. Primary 60G99, secondary
60G22, 60F17 .} \footnotemark[0]\footnotetext[0]{ \textit{Key words
and phrases}. Random linear processes, limit theorems, tapered distributions }
\footnotemark[0]\footnotetext[0]{ \textit{Corresponding author:}
Vygantas Paulauskas, Department of Mathematics and
 Informatics, Vilnius university, Naugarduko 24, Vilnius 03225, Lithuania,
 e-mail:vygantas.paulauskas@mif.vu.lt}
}


\author{ Vygantas Paulauskas$^{\text{\small 1}}$ \\
{\small $^{\text{1}}$ Vilnius University, Department of Mathematics
and
 Informatics,}\\}



\maketitle

\begin{abstract}

In the paper we consider the partial sum process $\sum_{k=1}^{[nt]}X_k^{(n)}$, where $\{X_k^{(n)}=\sum_{j=0}^{\infty} a_{j}^{(n)}\xi_{k-j}(b(n)), \ k\in \bz\},\ n\ge 1,$ is a series of linear processes with tapered filter $a_{j}^{(n)}=a_{j}\ind{[0\le j\le \l(n)]}$ and  heavy-tailed tapered innovations $\xi_{j}(b(n), \ j\in \bz$. Both tapering parameters  $b(n)$ and $\l(n)$ grow to  $\infty$ as $n\to \infty$. The limit behavior of the partial sum process depends on the growth of these two tapering parameters and dependence properties of a linear process with non-tapered filter $a_i, \ i\ge 0$ and non-tapered innovations. We consider the case where $b(n)$ grows relatively slow (soft tapering), and all three cases of growth of $\l(n)$ (strong, weak, and moderate tapering). In these cases the limit processes (in the sense of convergence of finite dimensional distributions) are Gaussian.

\end{abstract}
\section{Introduction }

In \cite{Paul21} the limit theorems for sums of values of a linear process with tapered innovations were considered. Let
\begin{equation}\label{linpr1}
X=\{X_k, \ k\in \bz \}, \ X_k=\sum_{j=0}^\infty a_j\xi_{k-j},
\end{equation}
where the filter $\{a_j, j\ge 0\}$ and innovations $\xi_i, i\in \bz$, independent and identically distributed (i.i.d.) random variables, are such that the linear random process in (\ref{linpr1}) is well defined. Let us define
\begin{equation}\label{sum}
S_{n}=S_{n}(X):=\sum_{k=1}^{n} X_{k} \ {\rm and} \ S_{n}(t; X)=S_{[nt]}(X), \quad t\ge 0.
\end{equation}
 By means of linear processes, depending on moments of innovations, we can model stationary sequences with finite or infinite variance, while the properties of a filter of a linear process allow to model different dependence of a sequence $\{X_k\}$, namely, long-range, short-range, and negative dependence (we shall use the abbreviations LRD, SRD, and ND, respectively).
The asymptotic behavior of partial sum processes $S_{n}(t; X)$ is well investigated,  starting with pioneering works  of Davydov \cite{dav1970} (the case of finite variance) and Astrauskas \cite{Astrauskas} (the case of infinite variance) and with a big list of subsequent papers, devoted to limit theorems for $S_{n}(t; X)$.

In many fields, especially connected with applications, we face random quantities distributed according power-law. For example, it is well-known that many  natural hazards,  such as earthquakes, rock falls, landslides, riverine floods, tsunami, wildfire exhibit power-law behavior, see, for example, \cite{Geist}.  Asserting that many processes in real life can be modeled using heavy-tailed distributions, at the same time we must admit that essentially real world is bounded, i.e., the quantities, which we are interested in, have some bounds or at least have much lighter tails. Good example is seismology. For many years the scalar value of seismic moment of an earthquake was modeled by Pareto distribution with some $x_0$ (minimal seismic moment, above which a seismic moment follows power-law) and exponent $\a$ (its density being $f(x)= \a x_0^\a x^{-\a-1},\ x\ge x_0 $). Theoretical reasoning, based on the theory of branching processes, predicts that this exponent is universal (for many types of earthquakes - ground, oceanic, shallow, etc.) constant equal $1/2$, although estimations from real data usually give a little bit bigger value, see discussion in \cite{Kagan}. Another important fact was that empirical data of earthquakes showed that while seismic moment follows power-law in quite big range of values, the largest values of samples demonstrate much lighter tails, see Figure 1 in \cite{Kagan}. Therefore it was suggested to apply the exponential taper to Pareto distribution.  Strict definitions of tapered random variables and distributions will be defined in Section 3, here we provide only the main meaning of this notion. If $\xi$ has a heavy-tailed distribution with the tail index $\a$ and $R$ is a non-negative random variable, independent from $\xi$ and having the light tail, then tapered random variable $\xi(b)$ with the taper $R$ and tapering parameter $b>0$ is defined as follows
$$
\xi(b)=\xi \ind{[|\xi|<b]}+\frac{\xi}{|\xi|}(b+R)\ind{[|\xi|\ge b]}.
$$
Usually $R$ is exponential law, but it can be degenerated at zero random variable, then $\xi (b)$ will be truncated (at level $b>0$) random variable. An example, illustrating the last possibility can  be provided from computer science. In many problems it is assumed that the file size in networks has heavy-tailed distribution, on the other hand there is so-called the File Allocation Table (FAT, see Microsoft Knowledge Base Article 154997 (2007) or Wikipedia), used on most computer systems, which allows the
largest file size to be 4GB (more precisely, $2^{32}-1$ bytes). There are more examples where tapered or truncated (truncation can be considered as a particular case of tapering) heavy-tailed distributions are used to model various processes, we can recommend papers \cite{Aban} and \cite{Meerschaert} containing big lists of references with such examples.
Tapering parameter can be dependent on $n$, and if $b_n$ is growing to infinity comparatively slowly, we say that we have hard tapering, if it grows rapidly, then we have soft tapering (for strict definitions see Definition \ref{def1} in Section 3). The next natural step in investigation of $S_{n}(t; X)$ is, instead of a linear process $X$ with a fixed filter and a fixed sequence of innovations,  to consider the family of linear processes $X^{(n)}=\{X_k^{(n)}\}$, indexed by $n$ and for each $n$ taking as innovations a sequence of i.i.d.  copies of a tapered heavy-tailed random variable $\xi(b_n)$. Then one can expect that limit behavior of the appropriately normalized partial sum process $S_{n}(t; X^{(n)})$ depends on the sort of tapering: if we have hard tapering, then the limit will be Gaussian, while in the case of soft tapering stable processes will appear as limits. In the case where for each $n$ random variables $X_k^{(n)}, \ k\in \bz $ are i.i.d. (this will be if $a_0=1, a_j=0$ for $j\ge 1$) such answer is given in \cite{Chakrabarty}. Namely, in the case of hard tapering we have Gaussian limit law for appropriately normalized $S_{n}(1; X^{(n)})$, while in the case of soft tapering limit for $S_{n}(1; X^{(n)})$ will be the same as for sum $\sum_{k=1}^{n} \xi_{k}$ with the same centering and normalizing sequences. 
In papers \cite{Paul21}, \cite{Paul22}, and \cite{PaulDam5} limit theorems for linear random processes and linear random fields (with a specific structure of a filter) with tapered innovations were considered. The first attempt in \cite{Paul21}, where linear random processes were considered, was not very successful, the results were incomplete and there was even a mistake, see \cite{Paul23}. Then in \cite{Paul22} and \cite{PaulDam5}  considering limit theorems for linear
random fields with tapered innovations and changing some steps in the proofs, presented in \cite{Paul21}, in the case of linear random processes with tapered innovations we got   almost the final result, see Theorem 3 in \cite{PaulDam5}. We use the word "almost", since it remains  as open problem  what is a limit on the border between hard
and soft tapering. In the case of independent summands in \cite{Chakrabarty} this intermediate case was considered and it was shown that  some infinitely divisible law is a limit law in this case.

Linear processes are formed by the means of a filter and innovations. In \cite{Sabzikar} linear processes with tapered filters were introduced and limit theorems  for sums of values of such processes
were investigated.
In \cite{Paul21} it was mentioned that it is possible to consider not only tapered innovations, but together  tapered innovations and a filter, and in this paper we shall investigate this possibility. Now we introduce the notion of a linear process with a tapered filter, and our definition will be a little bit more general than definition in \cite{Sabzikar}. By a  taper for the filter of a linear process (\ref{linpr1}) we call a sequence $f(j)>0, \ j\ge 0,$ which is rapidly
decreasing as $j\to \infty$ and  $\sum_{j=0}^\infty |a_jf(j)|<\infty$.
 From the analogy with tapering of innovations it is natural to take sequence $f(j)$ exponentially decreasing. We use the (unusual) letter $f$ for a sequence for the  reason, that it is convenient to think that this sequence is obtained as values of some function, defined on $\bbr_+$, at integer points, for example $f(j)=\exp (-j), j\in \NN,$ and $f(x)=\exp (-x), x\in \bbr_+$. Then a linear process with a tapered filter is defined
$$
X_k=\sum_{j=0}^\infty a_jf(j)\xi_{k-j}, \ k\in \bz.
$$
Such linear process with a tapered filter, independently from initial filter $\{a_j\}$, is always with short-range dependence, due to the condition  $\sum_{j=0}^\infty |a_jf(j)|<\infty$, therefore, considering sums of values of a linear process, it is interesting to consider family of linear processes with tapers $f_n(j)$, depending on $n$. Additionally, it is natural to require that $f_n(j)$ is close to $1$ for $j=0, 1, \dots, k_n$ with some  $k_n \to \infty$, as
 $n\to \infty$.
Thus, we shall consider a family of linear processes
\begin{equation}\label{trunkfil}
X^{(n)}=\{X_k^{(n)} \}, \quad {\rm where} \ \  X_k^{(n)}=\sum_{j=0}^\infty a_jf_n(j)\xi_{k-j}, \ k\in \bz.
\end{equation}
In \cite{Sabzikar} the following family of tapers was used $f_n(j)=\exp (-j/\l(n))$ with $\l(n)\to \infty$ as $n\to \infty$. If there exists the limit
$$
\lim_{n\to \infty}\frac{n}{\l(n)}=\l_*,
$$
 then the family $X^{(n)}$ is {\it called strongly, weakly, or moderately} tapered, if $\l_*=\infty, \l_*=0,$ or $0<\l_*<\infty,$ respectively. Let us note that in \cite{Sabzikar} the term "tempered"was used instead of "tapered" and the taper was of the form $\exp (-j{\bar \l}(n))$ with ${\bar \l}(n)\to 0$ as $n\to \infty$. Since we shall deal with tapered innovations and filter together, for us it is more convenient to have both parameters of tapering unboundedly increasing. Also we  take a little bit different family of tapers, namely,
\begin{equation}\label{newtaper}
f_n(j)=\ind{[0\le j\le \l(n)]}.
\end{equation}
We think that this new taper (\ref{newtaper}) (which essentially is the truncation) is more natural, than that used in \cite{Sabzikar}. If tapering by exponential taper of innovations, which are random variables,  is quite natural (in many cases more natural than truncation), for filter coefficients (non random quantities) more natural is truncation, since in real life we always have only finite number of filter coefficients, only this finite number can be big. Although the main goal of the paper is to consider linear random processes with tapered filter and innovations, in Section 2 we consider the case where only the filter is tapered. This allow to see how results depend on the choice of the tapering family $f_n(j)$.
Thus, in Section 2 we consider limit theorems for family of linear processes (\ref{trunkfil}) with the taper from (\ref{newtaper}) but with non-tapered innovations with finite variance. We shall take
$$
\l(n)=cn^{\g_1} \quad {\rm with}\ \ 0<c<\infty, \ 0<\g_1<\infty.
$$
We have strongly, weakly, or moderately tapered family of linear processes (\ref{trunkfil}) (with taper from (\ref{newtaper})), if $0<\g_1<1, \ \g_1>1,$ or $\g_1=1$, respectively. It is clear that in the cases of strong and weak tapering the constant $c$ does not play any role, therefore in the cases $\g_1\ne 1$ we shall take $c=1$. In the case of moderate tapering  $\l_*=c^{-1}$, and this constant is important. Results, obtained in Section 2 are similar to those obtained in Section 4 in \cite{Sabzikar}. The only difference is in the case of moderate tapering, limit process in this case in Theorem \ref{thm1} is different from that in \cite{Sabzikar}.

Finally, in Section 3 we consider family of linear processes with tapered innovations and filter
\begin{equation}\label{trunkfilinn}
{\bar X}^{(n)}=\{{\bar X}_k^{(n)} \}, \quad   {\bar X}_k^{(n)}=\sum_{j=0}^\infty a_jf_n(j)\xi_{k-j}(b(n)), \ k\in \bz,
\end{equation}
where filter taper is from (\ref{newtaper}) and $\xi_{k}(b(n)), \ k\in \bz,$ are tapered innovations, defined in Section 3, see (\ref{paretodf}) and (\ref{tapparetodf}). We  take, as in \cite{Paul21},  $b(n)=n^\g$ and we  have hard, soft or intermediate tapering, if $0<\g<1/\a, \ \g>1/\a$, or $\g=1/\a$, respectively. We left the same terminology, used in \cite{Sabzikar} for tapering of the filter, since now we can say, for example, that we consider a linear process with weak and hard tapering, understanding that we have weak tapering of the filter and  hard tapering of innovations.
In the paper we consider only the case of hard tapering and all three types of filter tapering and three types of dependence, since in all these 9 cases limit processes are Gaussian and the proofs are based on calculation of variances. The case of soft tapering, where stable laws appear in limit and different techniques are used, is postponed for the subsequent paper (similar situation as with papers \cite{Paul22} and \cite{PaulDam5}).

\section{Limit theorems for linear processes with tapered filters}\label{sec1}

In this section we consider the family of linear processes
$$
X^{(n)}=\{X_k^{(n)} \}, \quad   X_k^{(n)}=\sum_{j=0}^\infty {\tilde a}_j\e_{k-j}, \ k\in \bz,
$$
where
$$
{\tilde a}_i= {\tilde a}_i^{(n)}=\left \{\begin{array}{ll}
              a_i, & {\rm if} \ 0\le i\le \l(n), \\
              0, & {\rm if} \ i> \l(n).
               \end{array} \right.
$$
 and $\e_{k}, \ k\in \bz,$ are i.i.d. random variables with $E\e_1=0, \ E\e_1^2=1.$
  Let us denote
$$
Z_n(t)=A_n^{-1} S_n(t, X^{(n)} ),
$$
where $A_n$ is a normalizing sequence for $ S_n(1, X^{(n)})$ and $S_n(t, X)$ is defined in (\ref{sum}). If the linear process $X_k$ is with non-tapered filter, that is $f_n(j)\equiv 1$ for each $n$ (or $\l(n)\equiv \infty$ in (\ref{newtaper}) for all $n$), then it is well-known that the limit process (in the sense of the convergence of finite-dimensional distributions (f.d.d.)) for $Z_n(t)$ is fractional Brownian motion (FBM) $B_H$ with Hurst parameter $0<H<1$, depending on some properties of the filter of $X_k$. Since in our work we do not seek results under the most general conditions on the filter, as in previous our papers \cite{Paul21}, \cite{Paul22}, and \cite{PaulDam5}, we assume that
\begin{equation}\label{cond1}
a_n\sim n^{-\b},
\end{equation}
where $\b >1/2$ (this condition ensures the correctness of definition of a linear process $X=\{X_k, \ k\in \bz\}$ with non-tapered filter; for the existence of $X_k^{(n)}$ for each $n$ this condition is irrelevant, later we shall discuss what is happening when $\b$ is even negative). Together with (\ref{cond1}) we shall consider three main sets of conditions, giving rise to three different types of dependence of the process $X$ and at the same time three types of memory, using classification proposed in \cite{Paul20}.

(i) $1/2<\b<1$ - the case of LRD and positive memory;

(ii) $\b>1$ and $\sum_{i=1}^\infty a_i<\infty$, the case of SRD and zero memory;

(iii) $1<\b<3/2$ and $\sum_{i=1}^\infty a_i=0$, the case of ND and negative memory.

Then it is well-known that
$$
\{n^{-H}S_n(t, X)\}\fdd \{B_H (t)\}
$$
where
$$
H= \left \{\begin{array}{ll}
                     \frac{3}{2}-\b  \ & {\rm in \ the \ cases} \ (i) \ {\rm and}\ (iii), \\
                      \frac{1}{2}    \ & {\rm in\ the\ case} \ (ii).
                                                                       \end{array} \right.
$$
Intuitively it is clear (and this is confirmed by results in \cite{Sabzikar}) that in the case of weak tapering we shall get the same FBM as in the case of untapered filter, while in the case of strong tapering in all cases of dependence, as a limit process we shall get Brownian motion (BM) $B(t):=B_{1/2}(t)$. Only in the case of moderate tapering we can get a limit process, different from FBM, as it is in \cite{Sabzikar}, where in the case of exponential taper  there was the  tempered fractional Brownian motion of the second kind as a limit process.

Since for a filter we have three cases of tapering and three cases of dependence, thus we  consider nine cases. We shall number them by index $j=1,\dots , 9$ in the following way: $j=1, 2, 3$ we attribute to strong tapering ($0<\g_1<1$) and three dependence cases (i), (ii), and (iii), respectively, indices $j=4, 5, 6$ are attributed to week tapering, and $j=7, 8, 9$ are given to the case of moderate tapering and three dependence types. For example, $j=8$ means that we consider moderate tapering of  a filter and SRD. Therefore, we shall use the notations $Z_n^{(j)}(t), (A_n^{(j)})^2, S_n^{(j)}(t, X^{(n)})$, but in cases where it is clear which index $j$ is considered, we shall skip this index from the notation. Let $\{U_n\}\fdd \{U_0\}$ stand for the convergence of processes $U_n$ to a process $U_0$ in the sense of the  f.d.d.
Our main result of this section is the following theorem.

\begin{teo}\label{thm1} Suppose that there exists $0<\d\le 1$ such that $E|\e_1|^{2+\d}<\infty$. Then, for all $j=1,\dots , 9$, we have
$$
\{Z_n^{(j)}(t)\}\fdd \{U^{(j)}(t)\},
$$
where normalizing sequences $(A_n^{(j)})$ and Gaussian limit processes $U^{(j)}(t)$ are  defined by means of their covariance functions, which  are given in Proposition \ref{prop1}. Particularly, $U^{(j)}(t)=B(t)$ for $j=1, 2, 3,5, 8,$ and $U^{(j)}(t)=B_H(t)$, \ $H=3/2-\b$,  for   $j=4, 6$. Processes $U^{(7)}(t)$ and $U^{(9)}(t)$ will be discussed at the end of this section.
\end{teo}
\begin{remark}\label{rem1} Most probably the statement of Theorem \ref{thm1} holds for $\d=0$, if in the proof we should use Lindeberg type condition. But it will require some calculations in all nine cases. On the other  hand, in \cite{Paul21} and \cite{Paul22}, dealing with tapered innovations (having moments of all orders) we used Lyapunov fractions of order $2+\d$. Thus, in order to use the results from these two cited above papers we additionally assumed the existence of the moment of the order $2+\d$ of innovations.
\end{remark}

{\it Proof of Theorem \ref{thm1}}
in order to prove limit theorem for $Z_n(t)$ we shall use the same scheme of the proof as in \cite{Paul21}: we calculate ${\rm Var}S_n(t, X^{(n)})$, then we get $A_n^2={\rm Var}S_n(1, X^{(n)})$,  find a limit $\lim_{n\to \infty}{\rm Var}Z_n(t)$, and, finally, prove that f.d.d. of $Z_n(t)$ are asymptotically Gaussian.

\begin{prop}\label{prop1} For all $j=1,\dots , 9$ and $s, t >0$ we have
\begin{equation}\label{VarZnfinal}
\lim_{n\to \infty}{\rm Var}Z_n^{(j)}(t)= W^{(j)}(t),
\end{equation}
\begin{equation}\label{CovZnfinal}
\lim_{n\to \infty}{\rm Cov}\left (Z_n^{(j)}(t), Z_n^{(j)}(s)\right )= W^{(j)}(t)+W^{(j)}(s)-W^{(j)}(|t-s|),
\end{equation}
where
$$
W^{(j)}(t)= t^{2H(j)}, \ {\rm for} \  j=1, 2, 3, 4, 5, 6, 8,
$$
$$
W^{(7)}(t)= t^{2H(7)}D_1(t, \b, c), \ \ W^{(9)}(t)= t^{2H(9)}D_2(t, \b, c), 
$$
Here
$$
H(j)=1/2,\ {\rm for} \  j=1, 2, 3,  5,  8, H(j)=3/2-\b, \ {\rm for} \  j=4, 6, 7, 9,
$$
 and $D_1(t, \b, c)=C_{19}(t, \b, c), \ D_2(t, \b, c)=C_{26}(t, \b, c)$. Normalizing sequences are defined as follows: $(A_n^{(j)})^2=\left (\sum_{i=0}^\infty a_{i} \right )^2 n$ for $j=2, 5, 8,$ \ $(A_n^{(j)})^2= C(\b)n^{1+2\g_1(1-\b)}$ for $j=1, 3$, \ $(A_n^{(j)})^2= C(\b)n^{2H}$ for $j=4, 6,$ and \ $(A_n^{(j)})^2= C(\b, c)n^{2H}$ for $j=7, 9.$
\end{prop}

{\it Proof of Proposition \ref{prop1}}
Using formulae (2.2) and (2.3) from \cite{Paul21} and recalling that $E\e_1^2=1$  we can write
$$
S_n(t, X^{(n)})=\sum_{j=-\infty}^{[nt]} d_{n,j,t}\e_{j},
$$
$$
{\rm Var}S_n(t, X^{(n)})=\sum_{j=-\infty}^{[nt]}| d_{n,j,t}|^2,
$$
where  $d_{n,j,t}=\sum_{k=1}^{[nt]}{\tilde a}_{k-j}$ for $j\le 0$ and $d_{n,j,t}=\sum_{k=j}^{[nt]}{\tilde a}_{k-j}$ for $j> 0$. As in \cite{Paul21} we can write
\begin{equation}\label{sumdnj}
\sum_{j=-\infty}^{[nt]} |d_{n,j,t}|^2=V_1(t) +V_2(t):=\sum_{j=-\infty}^0 |d_{n,j,t}|^2 + \sum_{j=1}^{[nt]} |d_{n,j,t}|^2.
\end{equation}
To find the asymptotic of $V_1(t), V_2(t)$ with respect to n, as in \cite{Paul21} instead of condition (\ref{cond1}), which means that $a_n=n^{-\b}(1+\d(n))$, where $\d(n)\to 0$, as $n\to \infty$, we can simply assume $a_n=n^{-\b}$ for $n\ge 1$. Although in the proof of this step some changes are needed, since in \cite{Paul21} the non-tapered filter was considered, but since $\l(n)\to \infty$ in all cases of tapering, there is no principal difficulties.

We start with the cases $0<\g_1<1$ (strong tapering) and (i) (LRD), i.e. with the case $j=1$ (but, as it was mentioned above, we skip the index $j$ from notation). In the same way as in \cite{Paul21}, changing sums into integrals and denoting $m=[nt], m_1=[n^{\g_1}]$, we can get
$$
V_1(t)=\sum_{j=0}^\infty \left (\sum_{k=1}^m a_{k+j}\ind{[0\le k+j\le m_1]} \right )^2=\sum_{j=0}^{m_1} \left (\sum_{k=1}^{m_1-j} a_{k+j}\right )^2
$$
and
$$
V_1(t)\sim \int_0^{m_1}\left (\int_0^{m_1-y}(x+y)^{-\b} dx\right )^2 dy=\frac{1}{(1-\b)^2}m_1^{3-2\b}\int_0^1(1-v^{1-\b})^2)dv.
$$
Here we used the fact that for any fixed $t>0$ for sufficiently large $n$ we have $m_1<m$. For the second term in (\ref{sumdnj}) we have
$$
V_2(t)=\sum_{j=1}^m \left (\sum_{i=0}^{m-j} a_{i}\ind{[0\le i\le m_1]} \right )^2=V_{2,1}(t)+V_{2,2}(t),
$$
where
$$
V_{2,1}(t)=\sum_{j=1}^{m-m_1} \left (\sum_{i=0}^{m_1} a_{i} \right )^2, \quad V_{2,2}(t)=\sum_{j=m-m_1+1}^{m} \left (\sum_{i=0}^{m-j} a_{i} \right )^2.
$$
Similarly, as in the estimation of $V_1(t)$, we get
$$
V_{2,1}(t)\sim \int_0^{m-m_1}\left (\int_0^{m_1}x^{-\b} dx\right )^2 dy=\frac{(m-m_1)m_1^{2(1-\b)}}{(1-\b)^2},
$$
$$
V_{2,2}(t)\sim \int_{m-m_1}^m\left (\int_0^{m-y}x^{-\b} dx\right )^2 dy=\frac{m_1^{3-2\b}}{(1-\b)^2(3-2\b)}.
$$
From all these relations we get
\begin{equation}\label{asV1V2t}
V_1(t)+V_2(t)\sim \frac{(m-m_1)m_1^{2(1-\b)}}{(1-\b)^2}+\frac{m_1^{3-2\b}}{(1-\b)^2(3-2\b)}.
\end{equation}
The main term from these two terms is the first one, therefore we  take $t=1$ in (\ref{asV1V2t}) and $A_n^2=(1-\b)^{-2}n^{1+2\g_1(1-\b)}$ (we recall that $m_1/m\to 0$, as $n\to \infty$). Then we get
\begin{equation}\label{VarZn}
\lim_{n\to \infty}{\rm Var}Z_n(t)=t.
\end{equation}

In the case $j=3$, that is in the case  of negative dependence (iii), in estimation of $V_1(t)$ we use the fact that $\int_0^1(1-v^{1-\b})^2)dv<\infty$, since $1<\b<3/2$, while in estimation of $V_2(t)$ we use the relation $\sum_{i=0}^{k} a_{i}=-\sum_{i=k+1}^\infty a_{i}$, and in a similar way, as in the case (i), we can get  the same relation (\ref{asV1V2t}). Then, taking the same $A_n^2$, we  get the relation (\ref{VarZn}).

It remains in the case of strong tapering the case of short-range dependence (ii), the case $j=2$. Again, the main term is $V_{2,1}(t)$, and we easily get
$$
\frac{V_{2,1}(t)}{n}=\frac{1}{n}\sum_{j=1}^{m-m_1} \left (\sum_{i=0}^{m_1} a_{i} \right )^2=\frac{m-m_1}{n}\left (\sum_{i=0}^{m_1} a_{i} \right )^2 \to t\left (\sum_{i=0}^\infty a_{i} \right )^2,
$$
as $n\to \infty$. It is easy to get the following relations
$$
\frac{V_{2,2}(t)}{n}=\sum_{j=m-m_1}^m \left (\sum_{i=0}^{m-j} a_{i} \right )^2\le \frac{m_1}{n}\left (\sum_{i=0}^\infty |a_{i}| \right )^2 \to 0,
$$
$$
\frac{V_{1}(t)}{n}=\sum_{j=0}^{m_1} \left (\sum_{i=1}^{m_1-j} a_{i} \right )^2\le \frac{m_1}{n}\left (\sum_{i=0}^\infty |a_{i}| \right )^2 \to 0.
$$
Therefore, now we take $A_n^2=n\left (\sum_{i=0}^\infty a_{i} \right )^2$ and again we get (\ref{VarZn}).

Now we consider the case $\g_1 >1$ (weak tapering).  For any fixed $t>0$ and for sufficiently large $n$ we have the reverse inequality $m_1>m$ and $m_1/m\to \infty$. We start with the case (i), that is, $j=4$. Taking into account that $m_1>m$ we can write
$$
V_1(t)=\sum_{j=0}^\infty \left (\sum_{k=1}^{m} a_{k+j}\ind{[0\le k+j\le m_1]} \right )^2=V_{1,1}(t)+V_{1,2}(t),
$$
where
$$
V_{1,1}(t)=\sum_{j=0}^{m_1-m} \left (\sum_{i=1}^{m} a_{i+j} \right )^2, \quad V_{1,2}(t)=\sum_{j=m_1-m+1}^{m_1-1} \left (\sum_{i=0}^{m_1-j} a_{i+j} \right )^2.
$$
It is not difficult to get
\begin{eqnarray*}
\frac{V_{1,1}(t)}{n^{3-2\b}} &\sim& \frac{1}{n^{3-2\b}}\int_0^{m_1-m}\left (\int_0^{m}(x+y)^{-\b} dx\right )^2 dy \\
          &=&\left (\frac{m}{n}\right )^{3-2\b} \int_0^{(m_1-m)/m}\left (\int_0^{1}(x+y)^{-\b} dx\right )^2 dy .
\end{eqnarray*}
Since $C_1(\b=)\int_0^{\infty}\left (\int_0^{1}(x+y)^{-\b} dx\right )^2 dy<\infty$ for $1/2<\b<1$, therefore we have
\begin{equation}\label{asV11t}
\frac{V_{1,1}(t)}{n^{3-2\b}} \to C_1(\b)t^{3-2\b}.
\end{equation}
Similarly we get
$$
\frac{V_{1,2}(t)}{n^{3-2\b}}\sim \left (\frac{m}{n}\right )^{3-2\b}\int_{(m_1-m+1)/m}^{(m_1-1)/m}\left (\int_0^{(m_1/m)-y}(x+y)^{-\b} dx\right )^2 dy
$$
Let us denote $B_n=m_1/m$, then we shall show that
\begin{equation}\label{Dnb}
D(n, \b)=\int_{B_n-1+1/m}^{B_n-1/m}\left (\int_0^{B_n-y}(x+y)^{-\b} dx\right )^2 dy \to 0 \ {\rm as} \ n\to \infty.
\end{equation}
Since the length  of the interval for the outer integral is less than $1$, we can estimate
\begin{eqnarray*}
D(n, \b) &<& \sup_{B_n-1+1/m\le y \le B_n-1/m}\left (\int_0^{B_n-y}(x+y)^{-\b} dx\right )^2 \\
     &\le & \left (\int_0^{1}(x+B_n-1+1/m)^{-\b} dx\right )^2 \le (B_n-1)^{-2\b} \to 0 \ {\rm as} \ n\to \infty.
\end{eqnarray*}
Therefore, we have
\begin{equation}\label{asV12t}
\frac{V_{1,2}(t)}{n^{3-2\b}} \to 0.
\end{equation}
It remains to investigate $V_2(t)$. We have
$$
V_2(t)=\sum_{j=1}^m \left (\sum_{i=0}^{m-j}  a_{i}\ind{[0\le i\le m_1]} \right )^2=\sum_{j=1}^m \left (\sum_{i=0}^{m-j}  a_{i} \right )^2=\sum_{k=0}^{m-1} \left (\sum_{i=0}^{k}  a_{i} \right )^2,
$$
then it is easy to get
\begin{equation}\label{asV2t}
\frac{V_{2}(t)}{n^{3-2\b}} \to \frac{t^{3-2\b}}{(1-\b)^2(3-2\b)}.
\end{equation}
Collecting (\ref{asV11t})-(\ref{asV2t}) and taking $A_n^2=C_2(\b)n^{3-2\b}$,  we get
\begin{equation}\label{VarZn2}
\lim_{n\to \infty}{\rm Var}Z_n(t)=t^{2H}, \ H=\frac{3}{2}-\b.
\end{equation}
Here
$$
C_2(\b)=\frac{1}{(1-\b)^2(3-2\b)}+C_1(\b).
$$

As in the case of strong tapering, the case (iii) is similar to the case (i), so in the case $j=6$ we  formulate the final result. Taking $A_n^2=C_3(\b)n^{3-2\b}$ with
$$
C_3(\b)=\int_0^1\left (\int_{1-v}^\infty u^{-\b}du \right)^2dv+C_1(\b)
$$
we again get  (\ref{VarZn2}), only we recall that now $1<\b<3/2$ and the finiteness of both integrals in the expression of $C_3(\b)$ was shown in \cite{Paul21}.

In the case $j=5$ ($\g_1>1$ and (ii)) it is easy to show that, taking $A_n^2=n\left (\sum_{i=0}^\infty a_{i} \right )^2$,   we get (\ref{VarZn}).

Now we consider the last and the most complicated case of moderate tapering, $\g_1=1$, since now $m_1=cn$ and $m_1\le m$ or $m_1\ge m$ if $c<t$ or $c>t$, respectively. We start with the case $j=8$, i.e. with the case (ii) (SRD), since this case is more simple comparing with (i) and (iii). Let us consider the case $0<t<c$, i.e., $m<m_1$. We have
$$
V_2(t)=\sum_{j=1}^m \left (\sum_{i=0}^{m-j} a_{i}\ind{[0\le i\le m_1]} \right )^2=\sum_{j=1}^m \left (\sum_{i=0}^{m-j} a_{i}\right )^2
$$
and
\begin{equation}\label{asV2A}
\frac{V_{2}(t)}{n} \to t\left (\sum_{j=0}^\infty a_j \right)^2.
\end{equation}
As in the case of the weak tapering, denoting ${\tilde n}=m_1-m$, we have
$$
V_1(t)=\sum_{j=0}^\infty \left (\sum_{k=1}^{m} a_{k+j}\ind{[0\le k+j\le m_1]} \right )^2=V_{1,1}(t)+V_{1,2}(t),
$$
where
$$
V_{1,1}(t)=\sum_{j=0}^{{\tilde n}} \left (\sum_{i=1}^{m} a_{i+j} \right )^2, \quad V_{1,2}(t)=\sum_{j={\tilde n}+1}^{m_1-1} \left (\sum_{i=1}^{m_1-j} a_{i+j} \right )^2.
$$
Taking into account the relation ${\tilde n}=n(c-t)(1+o(1))$, it is not difficult to  obtain the following estimate
$$
\frac{V_{1,1}(t)}{n} \le \frac{C (c-t)}{(1-\b)^2{\tilde n}}\sum_{j=0}^{{\tilde n}}\left ((j+1)^{2(1-\b)}+(j+m)^{2(1-\b)} \right ).
$$
From this estimate it is easy to infer that
\begin{equation}\label{asV11tA}
\frac{V_{1,1}(t)}{n} \to 0,
\end{equation}
since , if $\b>3/2,$ then $\sum_{j=0}^{\infty}(j+1)^{2(1-\b)}<\infty$, while if $1<\b<3/2$, then
$$
\frac{1}{{\tilde n}}\sum_{j=0}^{{\tilde n}}\left ((j+1)^{2(1-\b)}+(j+m)^{2(1-\b)} \right )\le \frac{{\tilde n}^{3-2\b} +(m+{\tilde n})^{3-2\b} -m^{3-2\b}}{(3-2\b){\tilde n}}
$$
and the right-hand side of this inequality tends to zero, since $2(1-\b)<0$. If $\b=3/2$ then there appears logarithmic function, but we have the same result.

In a similar way we can estimate
\begin{eqnarray*}
\frac{V_{1,2}(t)}{n} &=& \frac{1}{n}\sum_{j={\tilde n}+1}^{m_1-1} \left (\sum_{i=1}^{m_1-j} a_{i+j} \right )^2 \le \frac{1}{n}\sum_{j={\tilde n}+1}^{m_1}(m_1-j)^2 j^{-2\b} \\
              &\le & \frac{(m_1-{\tilde n}-1)^2}{n}\sum_{j={\tilde n}+1}^{m_1}j^{-2\b}\le \frac{C(\b)m^2{\tilde n}^{1-2\b}}{n},
\end{eqnarray*}
and again we get
\begin{equation}\label{asV12tA}
\frac{V_{1,2}(t)}{n} \to 0,
\end{equation}
Collecting (\ref{asV2A})-(\ref{asV12tA}) we have
\begin{equation}\label{asV1+V2}
\frac{V_{1}(t)+V_2(t)}{n} \to t\left (\sum_{j=0}^\infty a_j \right)^2.
\end{equation}
In the case $0<c\le t$, i.e., $m \le m_1$, in a similar way (only now denoting ${\tilde n}=m-m_1$) we get (\ref{asV1+V2}). Therefore, in the case $j=8$ taking $A_n^2=n\left (\sum_{j=0}^\infty a_j \right)^2 $ we get the same relation (\ref{VarZn}).

Let us consider the case of moderate tapering and LRD, the case $j=7$. We start assuming that $0<t\le c$, i.e., $m\le m_1$ (since $m=[nt], m_1=[nc]$ there can be equality even if $t<c$). In a similar way as in the case of weak tapering ($\g>1$), only taking into account that now $m_1/m \to c/t$, instead of $m_1/m \to \infty$, as $n\to \infty$, we  can get
\begin{equation}\label{asV1mod}
\frac{V_{1}(t)}{n^{3-2\b}} \to t^{3-2\b}\left (C_{4}(t, \b, c)+ C_{5}(t, \b, c) \right ),
\end{equation}
where
$$
C_{4}(t, \b, c)=\int_0^{(c/t)-1}\left (\int_0^{1}(x+y)^{-\b} dx\right )^2 dy,
$$
$$
C_{5}(t, \b, c)=\int_{(c/t)-1}^{(c/t)}\left (\int_0^{(c/t)-y}(x+y)^{-\b} dx\right )^2 dy.
$$
Similarly we get
\begin{equation}\label{asV2mod}
\frac{V_{2}(t)}{n^{3-2\b}} \to t^{3-2\b}C_{6}(\b), \ \ {\rm where} \ \ C_{6}(\b)=\int_0^1\left (\int_0^{1-y}x^{-\b} dx\right )^2 dy.
\end{equation}
From (\ref{asV1mod}) and  (\ref{asV2mod}) we have
\begin{equation}\label{asV1+V2mod}
\frac{V_1(t)+V_{2}(t)}{n^{3-2\b}} \to t^{3-2\b}C_{7}(t, \b, c),
\end{equation}
where
$$
C_{7}(t, \b, c)=C_{4}(t, \b, c)+ C_{5}(t, \b, c) + C_{6}(\b).
$$
In the case $0<c<t$, i.e., $m_1\le m$, in a similar way we can get
\begin{equation}\label{asV1+V2modA}
\frac{V_1(t)+V_{2}(t)}{n^{3-2\b}} \to t^{3-2\b}C_{11}(t, \b, c),
\end{equation}
where
$$
C_{11}(t, \b, c)=C_{8}(t, \b, c)+C_{9}(t, \b, c)+C_{10}(t, \b, c),
$$
$$
C_{8}(t, \b, c)=\int_0^{1-(c/t)}\left (\int_0^{(c/t)}x^{-\b} dx\right )^2 dy, \ C_{9}(t, \b, c)=\int_{1-(c/t)}^{1}\left (\int_0^{1-y}x^{-\b} dx\right )^2 dy,
$$
$$
C_{10}(t, \b, c)=\int_0^{(c/t)}\left (\int_0^{(c/t)-y}(x+y)^{-\b} dx\right )^2 dy.
$$
Combining (\ref{asV1+V2mod}) and (\ref{asV1+V2modA}) we have
\begin{equation}\label{asV1+V2modB}
\frac{V_1(t)+V_{2}(t)}{n^{3-2\b}} \to t^{3-2\b}C_{12}(t, \b, c),
\end{equation}
where
$$
C_{12}(t, \b, c)=\left \{\begin{array}{ll}
 C_{7}(t, \b, c) ,               & 0<t\le c, \\
 C_{11}(t, \b, c) , & \ 0<c<t.
  \end{array} \right.
$$
There is continuity of $C_{12}$ at point $t=c$, since $C_{4}(c, \b, c)=C_{8}(c, \b, c)=0$ and it is easy to verify that $ C_{5}(c, \b, c) + C_{6}(\b)=C_{9}(c, \b, c)+C_{10}(c, \b, c)$.
Now we can choose the normalizing sequence $A_n^2$, to this aim we put $t=1$ in (\ref{asV1+V2modB}). It is easy to see that
$$
C_{12}(1, \b, c)=\left \{\begin{array}{ll}
 C_{7}(1, \b, c) ,               & 1\le c, \\
 C_{11}(1, \b, c) , & \ 0<c<1.
  \end{array} \right.
$$
Therefore, taking $A_n^2= C_{12}(1, \b, c)n^{3-2\b}$ we get
\begin{equation}\label{VarZn3}
\lim_{n\to \infty}{\rm Var}Z_n(t)=t^{3-2\b}C_{13}(t, \b, c), \ {\rm where} \ \  C_{13}(t, \b, c)=\frac{C_{12}(t, \b, c)}{C_{12}(1, \b, c)}.
\end{equation}

It remains the case $j=9$, the case of moderate tapering and (iii). Since the case (iii) is similar to (i), only now we use the equality $\sum_{i=0}^{k} a_{i}=-\sum_{i=k+1}^\infty a_{i}$, we provide only the final result. In the case $0<t<c,  m\le m_1, $ we get
\begin{equation}\label{asV1+V2modiii}
\frac{V_1(t)+V_{2}(t)}{n^{3-2\b}} \to t^{3-2\b}C_{14}(t, \b, c),
\end{equation}
where
$$
C_{14}(t, \b, c)=C_{4}(t, \b, c)+ C_{5}(t, \b, c) + C_{15}(\b).
$$
In the case $0<c<t,  m_1\le m, $ we get
\begin{equation}\label{asV1+V2modAiii}
\frac{V_1(t)+V_{2}(t)}{n^{3-2\b}} \to t^{3-2\b}C_{16}(t, \b, c),
\end{equation}
where
$$
C_{16}(t, \b, c)=C_{10}(t, \b, c)+C_{17}(t, \b, c)+C_{18}(t, \b, c),
$$
The expressions of these constants are as follows:
$$
 C_{15}(\b)=\int_0^1\left (\int_{1-y}^\infty x^{-\b} dx\right )^2 dy   , \  C_{17}(\b)=\int_{1-(c/t)}^{1}\left (\int_{1-y}^\infty x^{-\b} dx\right )^2 dy,
$$
$$
C_{18}(t, \b, c)=\int_0^{1-(c/t)}\left (\int_{(c/t)}^\infty x^{-\b} dx\right )^2 dy.
$$
Combining (\ref{asV1+V2modiii}) and (\ref{asV1+V2modAiii}) we have
$$
\frac{V_1(t)+V_{2}(t)}{n^{3-2\b}} \to t^{3-2\b}C_{19}(t, \b, c),
$$
and, taking $A_n^2= C_{19}(1, \b, c)n^{3-2\b}$, we get
$$
\lim_{n\to \infty}{\rm Var}Z_n(t)=t^{3-2\b}C_{20}(t, \b, c), \ {\rm where} \ \  C_{20}(t, \b, c)=\frac{C_{19}(t, \b, c)}{C_{19}(1, \b, c)}.
$$
Here
$$
C_{19}(t, \b, c)=\left \{\begin{array}{ll}
 C_{14}(t, \b, c) ,               & 0<t\le c, \\
 C_{16}(t, \b, c) , & \ 0<c<t.
  \end{array} \right.
$$

Thus, for all cases $j=1, \dots , 9$ we had proved (\ref{VarZnfinal}), then (\ref{CovZnfinal}) easily follows. Proposition \ref{prop1} is proved.

The second step in the proof of the theorem  is to prove that f.d.d. of $Z_n^{(j)}(t)$ are asymptotically normal. As it was noted in Remark \ref{rem1}, we use Lyapunov fractions of order $2+\d$ and results from \cite{Paul21} and \cite{Paul22}.   The following proposition gives us the asymptotic normality of f.d.d. of $Z_n^{(j)}(t)$.
We recall that $m=[nt]$.
\begin{prop}\label{prop2} If $E\e_1=0, E|\e_1|^{2}=1$ and there exists $0<\d\le 1$ such that $E|\e_1|^{2+\d}<\infty$, then, for all $j=1, 2, \dots, 9,$ as $n\to \infty$,
\begin{equation}\label{Liapunov}
L^{(j)}(2+\d,n, t):=\frac{\sum_{k=-\infty}^m |d_{n,k,t}^{(j)}|^{2+\d}E|\e_1|^{2+\d}}{\left (\sum_{k=-\infty}^n (d_{n,k,1}^{(j)})^2\right )^{(2+\d)/2}} \to 0.
\end{equation}
\end{prop}
{\it Proof of Proposition \ref{prop2}}. From the proof of Proposition \ref{prop1}, for any fixed $t$ and for all $j=1, 2, \dots, 9,$ we can estimate
\begin{equation}\label{LiapunovA}
\frac{\sum_{k=-\infty}^m (d_{n,k,t}^{(j)})^2}{\sum_{k=-\infty}^n (d_{n,k,1}^{(j)})^2}\le C(t),
\end{equation}
then it is easy to see that (\ref{Liapunov}) will follow if
\begin{equation}\label{Liapunov1}
{\tilde L}^{(j)}(\d,n, t)=\left (\frac{\max_{-\infty<k\le m} |d_{n,k,t}^{(j)}|}{A_n^{(j)}}\right )^{\d} \to 0.
\end{equation}
This quantity with $\d=1$ was estimated in the case of non-tapered filter in \cite{Paul21}, the estimation in the case of tapered filter goes along the same lines as in \cite{Paul21}, therefore we shall provide the proof of (\ref{Liapunov1}) only for several $j$. Let us denote
\begin{equation}\label{Liapunov1A}
I_1^{(j)}=\max_{0\le k <\infty}\big |\sum_{i=1}^{m}{\tilde a}_{i+k}\big |, \quad I_2^{(j)}=\max_{0<k\le m}\big |\sum_{i=k}^{m}{\tilde a}_{i-k}\big |,
\end{equation}
then (\ref{Liapunov1}) will follow if we show
\begin{equation}\label{Liapunov2}
\left (A_n^{(j)} \right )^{-1}\max (I_1^{(j)}, I_2^{(j)}) \to 0.
\end{equation}
Let us consider the case $0<\g_1<1$, then, for any fixed $t$, $m_1=[n^{\g_1}]< m$ and we have, for $j=1, 2, 3$,
$$
I_1^{(j)}=\max_{0\le k <\infty}\Big |\sum_{i=1}^{m}a_{i+k} \ind{\{0\le i+k \leq m_1\}}\Big |=\max_{0\le k \le m_1-1}\Big |\sum_{i=1}^{ m_1-k}a_{i+k}\Big |\le \sum_{i=1}^{m_1}|a_{i}|
$$
In a similar way we get
$$
I_2{(j)}=\max_{0\le k \le m}\big |\sum_{i=k}^{m}{\tilde a}_{i-k}\big |=\max_{0\le k \le m}|\sum_{i=0}^{\min(m, m_1-k)}a_{i}|\le \sum_{i=0}^{ m_1}|a_{i}|.
$$
From these two estimates we get
$$
\max (I_1^{(1)}, I_2^{(1)})\le C m_1^{1-\b} , \ \max (I_1^{(j)}, I_2^{(j)})\le C, \ {\rm for} \ j=2, 3.
$$
Taking into account the expressions of $A_n^{(j)} $ from Proposition \ref{prop1} we get (\ref{Liapunov2})for $j=1, 2, 3.$

In the case $\g_1>1$ we have $m_1>m$ and in a similar way we can get the following estimates
$$
I_1^{(j)}=\max_{0\le k \le m_1-1}\Big |\sum_{i=1}^{\min (m, m_1-k)}a_{i+k}\Big |\le \max \left (\sum_{i=1}^{ m}|a_{i}|, \sum_{i=m_1-m}^{ m_1}|a_{i}| \right )
$$
and
$$
I_2{(j)}\le \max_{0\le k \le m}\Big |\sum_{i=0}^{\min (m_1, m-k)}a_{i}\Big |\le \sum_{i=0}^{ m}|a_{i}|.
$$
Again, from these two estimations and the expressions of $A_n^{(j)}$ we get (\ref{Liapunov2})for $j=4, 5, 6.$

In the case $\g_1=1$ and $m_1=[nc]$ we must consider two cases $m<m_1$ and $m\ge  m_1$ and in the same way as above we  get (\ref{Liapunov2}) for $j=7, 8, 9.$ Thus, we had proved (\ref{Liapunov}).
\halmos

Propositions \ref{prop1} and \ref{prop2} prove Theorem \ref{thm1}.
\halmos

As it was expected and as it was demonstrated in \cite{Sabzikar}, in most cases we have the same limit processes $B_H, \ 0<H<1,$ as in the case of non-tapered filter. Namely, in the cases $j=1, 2, 3,5, 8$ we have $U^{(j)}(t)=B(t)$ and  for $j=4, 6$ we have $U^{(j)}(t)=B_H(t)$, where $B_H$ is FBM  with the Hurst parameter $H$ and $B$ is BM. Only in the case of moderate tapering and cases LRD and ND as a limit process we get a Gaussian process, different from $B_H$. In \cite{Sabzikar1} there was introduced the following process, named as tempered fractional Brownian motion (TFBM),
\begin{equation}\label{TFBM}
B_{\a, \l}(t)=\int_{-\infty}^{\infty} g_{\a, \l, t}(x) M_B(dx),
\end{equation}
where $\a<1/2, \ \l \ge 0$,  $M_B$ is an independently scattered Gaussian random measure with Lebesgue control  measure, and
$$
g_{\a, \l, t}(x)=(t-x)_+^{-\a}\exp (-\l (t-x)_+)- ((-x)_+)^{-\a}\exp (-\l (-x)_+).
$$
In the case $\l=0, \ -\a=H-1/2$ the process $B_{\a, \l}(t)$ coincides with FBM $B_H(t)$. This process was introduced without any relation with linear processes with tapered filters. In \cite{Sabzikar} and \cite{Sabzikar2} in relation with limit theorems for a linear process with the tapered filter similar process was  introduced and named as tempered fractional Brownian motion of the second kind (as partial case of more general tempered fractional stable motion of the second kind for $1<\a\le 2$ and abbreviated as TFBMII):
\begin{equation}\label{TFBMII}
B_{H, \l}^{II}(t)=\int_{-\infty}^{\infty} h_{H, \l} (t; x) M_B(dx),
\end{equation}
where
\begin{eqnarray*}
h_{H, \l} (t; x) &=& (t-x)_+^{H-\frac{1}{2}}\exp (-\l (t-x)_+)- ((-x)_+)^{H-\frac{1}{2}}\exp (-\l (-x)_+) \\
   & + & \int_0^t (s-x)_+^{H-\frac{1}{2}}\exp (-\l (s-x)_+) ds.
\end{eqnarray*}
Variances of the processes defined in (\ref{TFBM}) and (\ref{TFBMII}), as functions of $t$, are provided in \cite{Sabzikar1} and \cite{Sabzikar2}, respectively. The expressions of these variances are rather complicated, in the first case involving the modified Bessel function of the second kind and generalized hypergeometric function in the second case, we do not provide them here. But it is possible  to get the asymptotic of these variances as $t\to 0$ or $t\to \infty$ and then to compare  with the asymptotic of functions $W^{(7)}(t)$ and $W^{(9)}(t)$ from (\ref{VarZnfinal}). Let us note that functions $D_1$ and $D_2$ in expressions of $W^{(7)}$ and $W^{(9)}$ in Proposition \ref{prop1} are the functions of $z=c/t$. The integrands in these expressions are  power functions, so all integrals, except one, can be integrated and we get sums of powers of $z$ and  $z-1$. Only integrating integral in the expression of $C_{4}$ we shall  get the integral $I(z)=\int_0^{z-1} (y(1+y))^{1-\b} dy$, for $z>1$. This integral, using tables of integrals, can be  written as
$$
I(z)=(z-1)^{2-\b}\sum_{k=0}^\infty \frac{(\b-1)_k}{2-\b+k}\frac{(-1)^k (z-1)^k}{k!},
$$
 where $(a)_k=a(a+1)\dots (a+k-1), k\ge 1, \ (a)_0=1$. But if we want to establish asymptotic of constant $C_{13}(t, \b, c)$ as $t\to 0$,
 it is better to consider asymptotic of $C_{7}(t, \b, c)$, as a function of $z$ and as $z\to \infty$. Looking for asymptotic of constant $C_{13}(t, \b, c)$,  as $t\to \infty$,  the   asymptotic of $C_{11}(t, \b, c)$, as $z\to 0$, must be considered. Let us consider $C_{7}(t, \b, c)$, and with some abuse of notation we shall write $C_j(z, \b)=C_j(t, \b, c)$ for $j=4, 5, 7, \dots, 12$. In \cite{Paul21}, see (2.6) and (2.7) therein, it was proved the finiteness of the integral $\int_0^{\infty}\left (\int_0^{1}(x+y)^{-\b} dx\right )^2 dy$ for $1/2<\b<3/2, \b\ne 1$, therefore we have
\begin{equation}\label{asC10}
C_{4}(z, \b) \to C(\b) \ \ {\rm as} \ z\to \infty.
\end{equation}
Here $C(\b)$ and in the sequel $C(\b, c)$ stand for  constants, depending only on parameters in parenthesis  and which can be different in different places, even in one formula. In a similar way we can show that
\begin{equation}\label{asC11}
C_{5}(z, \b)= (3-2\b)z^{2(1-\b)}+ O(z^{1-2\b}) \ \ {\rm as} \ z\to \infty.
\end{equation}
From (\ref{asC10}),  (\ref{asC11}), and the fact that $C_{6}$ depends only on $\b$, we get
\begin{equation}\label{asC13}
C_{7}(z, \b)= (3-2\b)z^{2(1-\b)}+C(\b)+ O(z^{1-2\b}) \ \ {\rm as} \ z\to \infty.
\end{equation}
Therefore, from (\ref{VarZn3}) and (\ref{asC13})  we derive that as $t \to 0$
\begin{equation}\label{VarZn5}
\lim_{n\to \infty}{\rm Var}Z_n^{(7)}(t)=C(\b, c)t +C(\b, c)t^{3-2\b}.
\end{equation}
In a similar way we can investigate $C_{11}(z, \b)$, as $z\to 0$, and to find that, as $t \to \infty$,
\begin{equation}\label{VarZn6}
\lim_{n\to \infty}{\rm Var}Z_n^{(7)}(t)=t\left (\frac{c^{2(1-\b)}}{(1-\b)^{2(1-\b)}} +C(\b)\frac{c}{t}\right ).
\end{equation}
 In a similar way it is possible in the case $j=9$ to investigate the asymptotic of $C_{19}(z, \b)$ as $z\to 0$ and $z\to \infty$ and to get the asymptotic behavior of $\lim_{n\to \infty}{\rm Var}Z_n^{(9)}(t)$ as $t\to 0$ and $t\to \infty$, similar to (\ref{VarZn5}) and (\ref{VarZn6}). From (\ref{VarZn5}) we see that variance of $U^{(7)}(t)$ behaves as $t$ for small $t$, while variances of both processes defined in (\ref{TFBM}) and (\ref{TFBMII}) behave as $t^H$. Similar situation is with $U^{(9)}(t)$, therefore
we can conclude that limit Gaussian processes in the cases $j=7, 9$ do not coincide with processes given in  (\ref{TFBM}) and (\ref{TFBMII}). On the other hand, the processes  $U^{(7)}(t)$ and $U^{(9)}(t)$ coincide with FBM in the case $c=\infty$ (in this case tapering function is identically equal to 1; parameter $c$ plays the same role as $\l$ plays in (\ref{TFBM}) and (\ref{TFBMII}), they are connected by the relation $c=\l^{-1}$), therefore these processes can be called tapered fractional Brownian motion of the third kind (TFBMIII). By the way, both words "tempered" and "tapered", which can be considered as synonyms give the same letter T in abbreviation. Both processes  $U^{(7)}(t)$ and $U^{(9)}(t)$ are defined by means of their variances and covariances, but it is possible to look for the representation of TFBMIII by means of stochastic integrals, and one may guess that  this representation would be similar to (\ref{TFBMII}), only integrand  must be modified  taking into account that now we use different tapering function (\ref{newtaper}). Also an interesting question is what kind of Gaussian process we can get in the same way modifying (changing exponential taper into indicator function) integrand function in (\ref{TFBM}). All these questions are left for the future research.

One more remark concerning Theorem \ref{thm1} is appropriate here. For a filter $\{a_i\}$ we assumed traditional condition $\b>1/2$ in (\ref{cond1}). But since we consider the tapered filter $\{{\tilde a}_i^{(n)}\}$ (essentially it is truncated filter;  we consider truncation as a particular case of tapering), it is possible to consider the case $\b<1/2$ and even $\b<0$, i.e. the coefficients of a filter unboundedly growing. Since we think that such filters are not realistic, we shall demonstrate only the simple case $\b=0$, which has some meaning: filter consists from finite, but big number of coefficients which are  all equal. Thus, we consider the linear random process ${\tilde X}^{(n)}=\{ {\tilde X}_{k}^{(n)}=\sum_{j=0}^{\infty} {\tilde a}_{j}^{(n)}\e_{k-j}, \ k\in \bz\}$ with the filter  ${\tilde a}_i^{(n)}=\ind{[0\le i\le \l(n)]}$,  $\l(n)=n^{\g_1}$, and three cases of tapering choosing parameter $\g_1$. We prescribe indexes $j=10, 11, 12$ to  cases of strong, weak, and moderate  tapering, respectively. Also we assume that innovations satisfy condition of Theorem \ref{thm1}. As earlier, denoting $m=[nt], m_1=[n^{\g_1}]$, we have
$$
{\rm Var}S_n(t, {\tilde X}^{(n)})=V_1(t)+ V_2(t),
$$
where
$$
V_1(t)=\sum_{j=0}^\infty \left (\sum_{k=1}^m \ind{[0\le k+j\le m_1]} \right )^2, \ V_2(t)=\sum_{j=1}^m \left (\sum_{i=0}^{m-j} \ind{[0\le i\le m_1]} \right )^2.
$$
Let  us consider the cases $0<\g_1<1$ (strong tapering, the case $j=10$), $m_1<m$. In this case it is not difficult to get
$$
V_1(t)\sim Cm_1^3, \quad V_2(t)\sim mm_1^2.
$$
This gives us the following relation
$$
{\rm Var}S_n^{(10)}(t, X^{(n)})\sim n^{1+2\g_1}t,
$$
therefore, taking $A_n^{(10)}=n^{\g_1+1/2}$, we get
$$
\lim_{n\to \infty}{\rm Var}Z_n^{(10)}(t)=t .
$$
In the cases $\g_1>1$ (weak tapering, the case $j=11$) in a similar way we get
$$
{\rm Var}S_n^{(11)}(t, X^{(n)})\sim n^{2+\g_1}t^2,
$$
and, taking $A_n^{(11)}=n^{1+\g_1/2}$, we get
$$
\lim_{n\to \infty}{\rm Var}Z_n^{(11)}(t)=t^2 .
$$
Little bit more complicated is the case $\g_1=1,  \ \l(n)=cn$ (moderate tapering, the case $j=12$), since now we must consider the cases $0<t<c$ and $0<c<t$, also $c<1$ and $c>1$. We skip all calculations and formulate the final result. Taking $\left (A_n^{(12)}\right )^2=C_{21}(c)n^{3}$, we get
\begin{equation}\label{VarZn11}
\lim_{n\to \infty}{\rm Var}Z_n^{(12)}(t)= \left \{\begin{array}{ll}
t^2 C_{22}(t,  c) ,               & 0<t\le c, \\
t C_{23}(t,  c) , & \ 0<c<t.
  \end{array} \right.
\end{equation}
Here $C_{21}(c)=c^2-c^3/3$, for $0<c\le 1$, and $C_{21}(c)=c-1/3$, for $c > 1$, and
$$
C_{22}(t,  c)=\frac{c(1-t/(3c))}{C_{21}(c)}, \ {\rm for} \ 0<t\le c,
$$
$$
C_{23}(t,  c)=\frac{c^2(1-c/(3t))}{C_{21}(c)}, \ {\rm for} \ 0<c< t.
$$
Having variances (and at the same covariances) of limit processes for $Z_n^{(j)}(t), j=10, 11, 12,$ it remains to prove the gaussianity of these limit processes. This can be done exactly as in Proposition  \ref{prop2}, by showing $L^{(j)}(2+\d,n, t) \to 0$ as $n\to \infty$ for all $j=10, 11, 12$. Thus we have proved the following
\begin{prop}\label{prop3}. Let ${\tilde X}_k^{(n)}, \ k\in \bz $ be a family of linear processes with the tapered filter $a_i^{(n)}=\ind{[0\le i\le \l(n)]}$, \ $\l(n)=n^{\g_1}$, and innovations satisfying the condition of Theorem \ref{thm1}. Then, for all $j=10, 11, 12,$, we have
$$
\{Z_n^{(j)}(t)\}\fdd \{U^{(j)}(t)\},
$$
where $U^{(10)}(t)=B(t),  U^{(11)}(t)=B(1)t$, and $ U^{(12)}$ is mean zero Gaussian process with variance, given in (\ref{VarZn11}).
\end{prop}

If in Theorem \ref{thm1} we considered values of $\b$ in the interval $1/2<\b <3/2, \b\ne 1$ (LRD and ND) and $\b>1$ (SRD), Proposition \ref{prop3} demonstrates that for tapered filters it is possible to consider the values $\b\le 1/2$. We have the following general picture. In the case of strong tapering ($0<\g_1<1$) we always get in limit the Brownian motion, despite how big is the growth of  coefficients of a filter (the case of negative $\b$; although in Proposition \ref{prop3} we considered only the case $\b=0,$ but it is easy to see that the same result we get for $\b<0$). This can be explained by the following fact. The variance of $X_k^{(n)}$ grows with $\b \to -\infty$ and with $\l(n)$, but this growth can be compensated by normalizing. Main factor in proving the asymptotic normality is dependence between summands in a sum $S_n(t, X^{(n)})$, and $X_k^{(n)}, k\ge 1$ for each $n$ are $\l(n)$-dependent and, if $\g_1<1$, the number of summands $nt$ grows more rapidly comparing with $\l(n)=n^{\g_1}$, therefore we  get the Brownian motion  as a limit. The same picture can be seen in \cite{Sabzikar}, Theorem 4.3 (i). Contrary, in the case of weak tapering ($\g_1>1$), the main factor becomes dependence between summands, since now $\l(n)$ grows more rapidly comparing with $nt$. In the case of LRD in Theorem \ref{thm1}, the case $j=4$ we have that $H=3/2-\b \to 1$, as $\b\to 1/2$, so one can expect that for all $\b\le 1/2$ the limit process will be degenerate $B(1)t$, as in the case $\b=0$ in Proposition \ref{prop3}. The case $\b<1/2$ is not considered in \cite{Sabzikar}, Theorem 4.3 (ii). More difficult to predict what is happening in the case of moderate tapering,  $\b<1/2$, and our tapering function (\ref{newtaper}). In the case of moderate tapering exponential tapering function is more convenient to work with, and in \cite{Sabzikar}, Theorem 4.3 (iii) for all values $\b<3/2$ the limit process is TFBMII from (\ref{TFBMII}) with parameter $H=3/2-\b>0$. One may expect that in the case of tapering function (\ref{newtaper}) and $\b<1/2$ the limit Gaussian process $U(t)$ will have variance behaving as $t^{2}$ for small $t$ and as $t$ for large $t$, behavior changing at the point $t=c$. Proposition \ref{prop3}, the case $j=12$ with $\b=0$ supports this expectation.

\bigskip

\section{Limit theorems for linear processes with tapered filters and innovations}\label{sec2}

In this section we investigate what happens if in all nine cases, considered in Theorem \ref{thm1}, we add assumption that innovations are heavy tailed and tapered. As it was mentioned in the Introduction, we  consider only the case of hard tapering with limit Gaussian processes. The case of soft tapering (with stable limit processes) is postponed for the subsequent paper, since this case requires different technique.
 Thus, we consider the family of linear processes with tapered innovations and filter, defined in (\ref{trunkfilinn}). As in Section \ref{sec1}, filter taper is from (\ref{newtaper}), and we assume the same condition (\ref{cond1}) and three types of dependence (SRD, LRD, ND), defined in Section \ref{sec1}. We shall consider the same nine cases  as in Section \ref{sec1}, only now instead of innovations with unit variance we shall consider heavy-tailed tapered innovations. We shall use tapered innovations which  were used in papers \cite{Paul21}, \cite{Paul22}, and \cite{PaulDam5}. Let $\theta =\theta (\a)$ stand for  the standard Pareto distribution with probability density and distribution functions
 \begin{equation}\label{paretodf}
f(x)= \a  x^{-\a-1}, \quad F(x)=1-\left (\frac{1}{x}\right )^\a, \ x\ge 1,
\end{equation}
where $\a>0$. Let  $R$ be the standard exponential random variable with the  density function $e^{-x}$, for $x\ge 0$, independent of $\theta$. Then  tapered (with tapering parameter $b>1$) standard Pareto random variable $\zeta= \zeta (\a, b)$  can be written as
\begin{equation}\label{tapparetodf}
\zeta (\a, b)=\theta\ind{[\theta<b]}+(b+R)\ind{[\theta\ge b]},
\end{equation}
its density function was written in \cite{Paul22}, see (1.4) therein.
Let us denote $\xi =\xi (\a, b)=\zeta (\a, b)-E\zeta (\a, b)$ . We consider the family (indexed by $n$) of linear random processes ${\bar X}^{(n)}=\{{\bar X}_{k}(b(n)),\ k\in \bz \}$,
\begin{equation}\label{tapinnfil}
{\bar X}_{k}(b(n))=\sum_{j=0}^{\infty} {\tilde a}_{j}^{(n)}\xi_{k-j}(b(n)), \ k\in \bz,
\end{equation}
where $\xi_{k}(b(n))$ are i.i.d. copies of a random variable $\xi (\a, b(n))$ with a tapering parameter $b(n)\to \infty$. (\ref{tapinnfil}) differs from (\ref{trunkfil})  only in change of innovations. All notations with bar  sign (${\bar Z}_n, {\bar A}_n,$ etc.) mean that ${\bar X}^{(n)}$ is used instead of  $X^{(n)}$. Here one remark is appropriate. Instead of standard Pareto random variable it is possible to taper more general random variable $\nu$, belonging to the domain of attraction of a stable random variable with exponent $\a$, as  it is done in \cite{Chakrabarty}, i.e., instead of $\zeta (\a, b)$ it is possible to consider the tapered random variable
$$
{\bar \zeta}(\a, b)=\nu\ind{[|\nu|<b]}+\frac{\nu}{|\nu|}(b+R)\ind{[|\nu|\ge b]}.
$$
But instead of considering this more general case, we, as in papers \cite{Paul21}, \cite{Paul22}, and \cite{PaulDam5}, continue with tapered Pareto distribution for the following reasons. Many calculations, considering $\nu$ instead of $\theta$, become more technical, we must require some asymptotic relations (with respect to $b\to \infty$) for $E|{\bar \zeta}(\a, b)|^r$, similar to (2.1) in \cite{Paul21}. Another reason is connected with applications. We know that the tail behavior of $\nu$ is $P\{|\nu|\ge x\}=x^{-\a}L(x)$, for large $x$, where $L$ is a slowly varying (at infinity) function,
while the tail of standard Pareto distribution is simply $x^{-\a}$, so the difference between these two tails is only slowly varying  function $L$. But from statistical theory we know that even the estimation of the tail index $\a$ is a very difficult problem. There are many estimators and hundreds of papers, devoted to this estimation, and practitioners face many difficulties when trying to estimate the tail index. And there is no hope from a sample to estimate what slowly varying  function $L$ is present in the tail of a distribution. If one looks at log-log plots of real data from various fields in \cite{Kagan}, \cite{Aban}, and \cite{Meerschaert} one can see the tapering effect, therefore during last decades there were papers dealing with estimation of parameters of Pareto tapered (including truncated) distributions, see the above cited papers \cite{Aban} and \cite{Meerschaert}. Tapered Pareto distribution, being the most simple between tapered other heavy-tailed distributions, quite well represents these other distributions, for example, in \cite{Meerschaert} there was generated  a sample from tapered  stable distribution and it turned out that tapered Pareto model gave quite good fit to this data. Of course, this fit can be explained by the fact that for stable laws the function $L(x)\to c$ as $x\to \infty$, it would be more interesting instead of  a stable law to take another heavy-tailed distribution with $L(x)\to \infty$ as $x\to \infty$. On the other hand, in \cite{Chakrabarty1} it is demonstrated what difficulties arise when trying to estimate parameters  of a tapered random variable $\nu$ even in the case of the function $L$ of the form $L(x)=C+O(x^{-\e})$ with some $\e>0$, i.e., $\nu$ belongs to the so-called Hall class.

Now we recall the notion of the hard and soft tapering, which was introduced in \cite{Chakrabarty} (only with the names of the hard and soft truncation) and was renamed and used in \cite{Paul21}.
\begin{definition}\label{def1} Let  $\{\theta_i \}, \ i\in N$ be i.i.d. random variables with regularly varying tails and with the tail exponent $0<a<2$. A sequence of tapering levels $\{b_n\}$ is called soft, hard, or intermediate tapering for a sequence $\{\theta_i, \ i\in N \}$,  if
$$
\lim_{n\to \infty} n\PP\{|\theta_1|>b_n\}
$$
 is equal to $0$,  $\infty,$ or $0<C<\infty$, respectively.
\end{definition}
 In our case the initial sequence $\{\theta_i\}, \ i\in N,$ is a sequence of standard Pareto random variables, therefore soft and hard tapering is defined as follows: if $b_n=n^\g$ and $\g>1/\a$, we have soft tapering, if $\g=1/\a$ -intermediate tapering, while $0<\g<1/\a$ gives us hard tapering. Clearly, only values $0<\a<2$ are interesting, and in the sequel we shall use this assumption without mentioning.

We consider
$$
{\bar Z}_n(t)={\bar A}_n^{-1} S_n(t, {\bar X}^{(n)}),
$$
where ${\bar X}^{(n)}$ is a family of linear processes with tapered innovations and filter, defined in (\ref{tapinnfil}). We shall consider the same nine cases, which were considered in Theorem \ref{thm1}, therefore in notations we add index $j$. The following result  shows that in the case of hard tapering of innovations the asymptotic behavior of ${\bar Z}_n^{(j)}(t)$ is the same as of $Z_n^{(j)}(t)$, given in Theorem \ref{thm1}.

\begin{teo}\label{thm2}  For all $j=1,\dots , 9$, we have
$$
\{{\bar Z}_n^{(j)}(t)\}\fdd \{U^{(j)}(t)\},
$$
where   Gaussian limit processes $U^{(j)}(t)$ are defined in Theorem \ref{thm1} and normalizing sequences ${\bar A}_n^{(j)}=A_n^{(j)}(E(\xi_{1}(b(n)))^2)^{1/2}$.
\end{teo}

{\it Proof of Theorem \ref{thm2}}.  As in the proof of Theorem \ref{thm1} there are two steps: calculation of ${\rm Var}S_n(t, {\bar X}^{(n)})$ and proving asymptotic normality of ${\bar Z}_n^{(j)}(t)$.
Since ${\rm Var}S_n(t, {\bar X}^{(n)})={\rm Var}S_n(t, X^{(n)})E(\xi_{1}(b(n)))^2$ and ${\bar A}_n^{2}=A_n^{2}E(\xi_{1}(b(n)))^2$, it is easy to see that from Proposition \ref{prop1} we get the following result.
\begin{prop}\label{prop4} For all $j=1,\dots , 9$ and $s, t >0$ we have
$$
\lim_{n\to \infty}{\rm Var}{\bar Z}_n^{(j)}(t)= W^{(j)}(t),
$$
$$
\lim_{n\to \infty}{\rm Cov}\left ({\bar Z}_n^{(j)}(t), {\bar Z}_n^{(j)}(s)\right )= W^{(j)}(t)+W^{(j)}(s)-W^{(j)}(|t-s|),
$$
where functions  $W^{(j)}$ are defined in Proposition \ref{prop1} and $({\bar A}_n^{(j)})^2=(A_n^{(j)})^2E(\xi_{1}(b(n)))^2$.
\end{prop}
In the previous section main difficulty in proving Theorem \ref{thm1} was calculation of ${\rm Var}S_n(t, X^{(n)})$, now the main difficulty is in proving asymptotic normality of ${\bar Z}_n^{(j)}(t)$, i.e., analog of Proposition \ref{prop2}. Instead of (\ref{Liapunov}) now we define
$$
{\hat L}^{(j)}(2+\d,n, t):=\frac{\sum_{k=-\infty}^m |d_{n,k,t}^{(j)}|^{2+\d}}{\left (\sum_{k=-\infty}^n (d_{n,k,1}^{(j)})^2\right )^{(2+\d)/2}}\frac{E|\xi_{1}(b(n))|^{2+\d}}{\left (E|\xi_{1}(b(n))|^{2}\right )^{(2+\d)/2}}.
$$
\begin{prop}\label{prop5} For all $j=1, \dots , 9$ there exists $0<\d=\d(j)\le 1$ such that, for $\g<1/\a$, as $n\to \infty$,
\begin{equation}\label{Liapunov3A}
{\hat L}^{(j)}(2+\d,n, t) \to 0.
\end{equation}
In the cases $j=1,2,4,5,7,8$ \ $\d(j)$ can be taken any number in the interval $(0, 1]$, while in the cases $j=3,6,9$ we must take  $0<\d(j)<\min (1, (3-2\b)(\b-1)^{-1})$.
\end{prop}
{\it Proof of Proposition \ref{prop5}}.
Taking into account formula (2.1) in \cite{Paul21} we have
\begin{equation}\label{Liapunov4}
\frac{E|\xi_{1}(b(n))|^{2+\d}}{\left (E|\xi_{1}(b(n))|^{2}\right )^{(2+\d)/2}}\le C(\a, \d)n^{\g \a \d/2}.
\end{equation}
Quantity $\sum_{k=-\infty}^n (d_{n,k,1}^{(j)})^2$ is $(A_n^{(j)})^2$ from Proposition \ref{prop1}, thus it remains to estimate the nominator of the  quantity
$$
{\tilde L}^{(j)}(2+\d,n, t)=\frac{\sum_{k=-\infty}^m |d_{n,k,t}^{(j)}|^{2+\d}}{\left (\sum_{k=-\infty}^n (d_{n,k,1}^{(j)})^2\right )^{(2+\d)/2}}.
$$
In some cases, namely, in the cases $j=2, 4, 5, 7, 8$, we can use the same method of estimation, which was used in Proposition \ref{prop2}:
\begin{equation}\label{Liapunov5}
{\tilde L}^{(j)}(2+\d,n, t)\le \left (\frac{\max_{-\infty<k\le m} |d_{n,k,t}^{(j)}|}{A_n^{(j)}}\right )^{\d}\frac{\sum_{k=-\infty}^m (d_{n,k,t}^{(j)})^2}{\sum_{k=-\infty}^n (d_{n,k,1}^{(j)})^2}.
\end{equation}
Using estimate (\ref{LiapunovA}) and notation (\ref{Liapunov1A}) we see that we need to show that
\begin{equation}\label{Liapunov6}
\left ( \left (A_n^{(j)} \right )^{-1}\max (I_1^{(j)}, I_2^{(j)})n^{\g \a /2}\right)^\d \to 0, \ {\rm as} \ n\to \infty.
\end{equation}
The quantity $\left (A_n^{(j)} \right )^{-1}\max (I_1^{(j)}, I_2^{(j)})$ was estimated in the proof of Proposition \ref{prop2}, therefore it remains to  take these estimates  and compare with $n^{\g \a /2}$. For example, in the case $j=2$ we have $\max (I_1^{(j)}, I_2^{(j)})\le C$ and $A_n^{(2)}$ is of the order $n^{1/2}$, therefore $-1/2 +\g \a /2<0$ if $\g<1/\a$. In a similar way we get (\ref{Liapunov6}) in the cases $j=4, 5, 7, 8$. Thus, we have proved (\ref{Liapunov3A}) for $j=2, 4, 5, 7, 8$, and in these cases $0<\d(j)\le 1$ can be any. But if we try the same method of estimation in the case $j=1,$ we shall get the following result
$$
 \left (A_n^{(j)} \right )^{-1}\max (I_1^{(j)}, I_2^{(j)})n^{\g \a /2}\le C(\a, \b)n^{(-1+2(1-\b)(1-\g_1)+\g\a)/2},
$$
and this quantity tends to zero if $\g\le (1-2(1-\b)(1-\g_1))/\a$. Since in this case $1/2<\b<1$ and $\g_1<1$,  we get more restrictive condition for $\g$. Therefore, we need to estimate ${\tilde L}^{(j)}(2+\d,n, t)$ without using (\ref{Liapunov5}). Let us consider the case $j=1$, i.e., the case $1/2<\b<1$ and $\g_1<1$, then
\begin{equation}\label{Jn12}
\sum_{k=-\infty}^m |d_{n,k,t}^{(1)}|^{2+\d}=J_1(n)+J_2(n),
\end{equation}
where (we skip in notations index $j=1$ and use the notations $m=[nt], m_1=[n^{\g_1}]$)
$$
J_1(n)=\sum_{i=0}^\infty \Big | \sum_{k=1}^m a_{k+i}\ind{[0\le k+i\le m_1]} \Big |^{2+\d}, \quad  J_2(n)=\sum_{i=1}^m \Big |\sum_{k=0}^{m-i} a_{k}\ind{[0\le k\le m_1]} \Big |^{2+\d}
$$
In Proposition \ref{prop1} these quantities were estimated in the case $\d=0$, on the other hand, in \cite{Paul22} such sums were estimated with $\d>0$ but with non-tapered filter. Therefore, in the same way we can get
\begin{eqnarray}\label{Jn1}
J_1(n) &=& \sum_{i=0}^\infty \Big | \sum_{k=1}^m a_{k+i}\ind{[0\le k+i\le m_1]} \Big |^{2+\d}=\sum_{i=0}^{m_1} \Big | \sum_{k=1}^{m_1-i} a_{k+i} \Big |^{2+\d} \nonumber \\
       & \le & \int_0^{m_1}\left (\int_0^{m_1-y}(x+y)^{-\b} dx\right )^{2+\d} dy \le C(\b, \d)m_1^{1+(1-\b)(2+\d)},
\end{eqnarray}
\begin{eqnarray}\label{Jn2}
J_2(n) &=& \sum_{i=1}^m \Big | \sum_{k=1}^{\min(m_1, m-i)} a_{k} \Big |^{2+\d}=\sum_{i=1}^{m-m_1} \Big | \sum_{k=1}^{m_1} a_{k} \Big |^{2+\d} +\sum_{i=m-m_1+1}^m \Big | \sum_{k=1}^{ m-i} a_{k} \Big |^{2+\d}  \nonumber \\
       & \le & \int_0^{m-m_1}\left (\int_0^{m_1}x^{-\b} dx\right )^{2+\d} dy + \int_{m-m_1}^m\left (\int_0^{m-y}x^{-\b} dx\right )^{2+\d} dy  \nonumber \\
       & \le & C(\b, \d, t)n^{1+\g_1(1-\b)(2+\d)}.
\end{eqnarray}
From (\ref{Jn12})-(\ref{Jn2}), taking into account that $\g_1<1$ we get
$$
{\hat L}^{(1)}(2+\d,n, t)\le C(\b, \d, t)n^{(-\d+\g\a\d)/2},
$$
and if $\g<1/\a$, then ${\hat L}^{(j)}(2+\d,n, t) \to 0$ as $n\to \infty.$ Thus we have proved (\ref{Liapunov3A}) in the case $j=1$, $\d(1)$ can be chosen arbitrary from the interval $(0, 1)$.
It remains three cases $j=3, 6, 9$, all with ND. Let us consider the case $j=3.$ As in the case $j=1$, we get
\begin{equation}\label{Jn31}
J_1(n)   \le m_1^{1+(1-\b)(2+\d)}\int_0^{1}\left (\int_0^{1-y}(x+y)^{-\b} dx\right )^{2+\d} dy,
\end{equation}
\begin{equation}\label{Jn32}
J_2(n)   \le m^{1+(1-\b)(2+\d)}\left ( C(\b, \d)(1+\int_{1-\frac{m_1}{m}}^{1}(1-y)^{(1-\b)(2+\d)}  dy)\right ).
\end{equation}
In the case $j=1$ we had $1/2<\b<1$, now we have $1<\b<3/2$, so we must show that integrals, appearing in (\ref{Jn31}) and (\ref{Jn32}),with appropriate choice of $\d$, are finite. As in \cite{Paul22} (see (2.20)-(2.22) therein) the integral in (\ref{Jn31}) is finite if $(1-\b)(2+\d)>-1$, or $0<\d<(3-2\b)(\b-1)^{-1}$. With the same bounds for $\d$ we get that the integral in (\ref{Jn32}) tends to zero, since $m_1m^{-1} \to 0$ as $n\to \infty.$ Therefore we can take
\begin{equation}\label{delta3}
\d=\min \left (1, \frac{3-2\b}{2(\b-1)}\right )
\end{equation}
and, collecting (\ref{Liapunov4}), (\ref{Jn12}), (\ref{Jn31}), and (\ref{Jn32}), we get for this value of $\d$
$$
{\hat L}^{(3)}(2+\d,n, t)\le C(\b, \d, t)n^{(-\d-2(1-\g_1)(\b-1)(2+\d) +\g\a\d)/2}.
$$
From this estimate we see that if
$$
\g<\frac{1}{\a}\left (1+\frac{2(1-\g_1)(\b-1)(2+\d)}{\d} \right )
$$
than we have (\ref{Liapunov3A}) in the case $j=3$, $\d(3)$ can be taken as in (\ref{delta3}).

Now we consider the case $j=6$, i.e., the case $\g_1>1$ and ND of the filter. Again, skipping the index $j=6$, we can get
\begin{eqnarray*}
J_1(n) &=& \sum_{i=0}^{m_1-m} \Big | \sum_{k=1}^{m} a_{k+i} \Big |^{2+\d} +\sum_{i=m_1-m+1}^{m_1-1} \Big | \sum_{k=0}^{m_1-i} a_{k+i} \Big |^{2+\d} \nonumber \\
       &\le & n^{1+(1-\b)(2+\d)} \int_0^{n^{\g_1-1}-t}\left (\int_0^{t}(x+y)^{-\b} dx\right )^{2+\d} dy \nonumber \\
       &+&  n^{1+(1-\b)(2+\d)} \int_{n^{\g_1-1}-t}^{n^{\g_1-1}} \left (\int_0^{n^{\g_1-1}-y}(x+y)^{-\b} dx\right )^{2+\d} dy   \nonumber \\
\end{eqnarray*}
The integral $ \int_0^{\infty}\left (\int_0^{t}(x+y)^{-\b} dx\right )^{2+\d} dy $ (only with $t=1$) was in the expression of $C_2(\b)$ in the previous section and is finite if $(1-\b)(2+\d)>-1$. Similarly to (\ref{Dnb}), we can prove that
$$
\int_{n^{\g_1-1}-t}^{n^{\g_1-1}} \left (\int_0^{n^{\g_1-1}-y}(x+y)^{-\b} dx\right )^{2+\d} dy \to 0 \  \ {\rm as} \ \ n\to \infty.
$$
Therefore, we get
\begin{equation}\label{Jn61}
J_1(n)   \le C(\b, \d, t)n^{1+(1-\b)(2+\d)}.
\end{equation}
Under the same condition $(1-\b)(2+\d)>-1$ we get
\begin{eqnarray*}
J_2(n) &=& \sum_{i=1}^{m} \big |\sum_{k=0}^{m-i} a_{k} \big |^{2+\d}=\sum_{i=1}^{m} \big |\sum_{k=m-i+1}^{\infty} a_{k} \big |^{2+\d} \nonumber \\
     &\le & n^{1+(1-\b)(2+\d)}\int_0^t\left (\int_{t-y}^\infty x^{-\b} dx\right )^{2+\d} dy. \nonumber \\
\end{eqnarray*}
Thus, we have
\begin{equation}\label{Jn62}
J_2(n)   \le C(\b, \d, t)n^{1+(1-\b)(2+\d)}.
\end{equation}
Now taking the value of $\d$ as in (\ref{delta3}) and collecting (\ref{Liapunov4}), (\ref{Jn61}), and (\ref{Jn62}),  we get for this value of $\d$
$$
{\hat L}^{(6)}(2+\d,n, t)\le C(\b, \d, t)n^{(-\d +\g\a\d)/2}.
$$
Thus, we have (\ref{Liapunov3A}) in the case $j=6$, $\d(6)$ can be taken as in (\ref{delta3}). The proof of the last case $j=9$ goes along the same lines as in the case $j=6,$ only now we must consider separately the cases $0<c\le t$ and $0<t<c$ and there will  appear the following integrals:
$$
\int_0^{c-t}\left (\int_0^{t}(x+y)^{-\b} dx\right )^{2+\d} dy, \ \  \int_{c-t}^{c} \left (\int_0^{c-y}(x+y)^{-\b} dx\right )^{2+\d} dy.
$$
Again, for the finiteness of these integrals we need the condition $(1-\b)(2+\d)>-1$, and we get (\ref{Liapunov3A}) in the case $j=9$, $\d(9)$ can be taken as in (\ref{delta3}). Proposition \ref{prop5} is proved.
\halmos

Propositions \ref{prop4} and \ref{prop5} prove Theorem \ref{thm2}.
\halmos


\end{document}